\documentclass[12pt]{amsart}
\usepackage{amsmath, amssymb, amsthm, amsfonts, mathrsfs, latexsym}
\usepackage{enumerate}

\makeatletter
\@namedef{subjclassname@2010}{%
  \textup{2010} Mathematics Subject Classification}
\makeatother

\newtheorem{thm}{Theorem}[section]
\newtheorem{thmIntro}{Theorem}
\newtheorem{cor}[thm]{Corollary}
\newtheorem{lemma}[thm]{Lemma}

\newtheorem*{qstIntro}{Question}
\newtheorem{fact}[thm]{Fact}

\theoremstyle{definition}
\newtheorem{defn}[thm]{Definition}
\newtheorem{rmk}[thm]{Remark}

\numberwithin{equation}{section}

\newcommand{\Ha}{\ensuremath{\mathbb{H} } }

\newcommand{\Aut}{\ensuremath{\mathrm{Aut} } }
\newcommand{\Cay}{\ensuremath{\mathrm{Cay} } }
\newcommand{\Th}{\ensuremath{\mathrm{Th} } }

\begin{document}


\baselineskip=17pt



\title[The random graph on Hall's universal group]
{The random graph on Hall's universal group}

\author[R.\,Sklinos]{Rizos Sklinos}

\author[S.\,Song]{Shichang Song}

\date{}

\begin{abstract}
For every integer $m\geq 7$, we show that the Cayley graph of Hall's universal locally finite group $\Ha$ with respect to the set of elements of order exactly $m$ is isomorphic to the random graph. In particular, the Fra\"{i}ss\'e limit of finite graphs is a parameter-free, quantifier-free definable reduct of the Fra\"{i}ss\'e limit of finite groups.

We use the order 7 reduct to code finite subsets of $\Ha$ uniformly and give an alternative proof of Hodges's theorem that standard arithmetic is interpretable without parameters in $\Ha$. Hall's construction gives a parameter-free interpretation in the reverse direction. We show that the two structures admit no bi-interpretation, even with parameters.

Finally, we show that the first order theory of $\Ha$ has no model companion and fails to geometrically eliminate imaginaries.
\end{abstract}



\maketitle

\section{Introduction}
Fra\"{i}ss\'e's construction is a central tool in model theory. A Fra\"{i}ss\'e class is assembled into a single countable structure that contains a copy of every member of the class. Every isomorphism between finite substructures extends to an automorphism of the entire limit, and the class determines this universal homogeneous structure up to isomorphism. The construction recovers structures that arise naturally in different areas of mathematics. Hall's universal locally finite group $\mathbb{H}$, the Rado graph $\Gamma$ (also called the random graph), and the rationals with their usual order are the limits of finite groups, finite graphs, and finite linear orders, respectively.

The first two examples suggest a question. 
\begin{qstIntro}
Can the random graph appear as a Cayley graph of Hall's universal group?
\end{qstIntro}

More explicitly, can we choose an inverse-closed generating set $S\subseteq\mathbb{H}\setminus\{1\}$ so that the graph with vertex set $\mathbb{H}$ and adjacency given by $x^{-1}y\in S$ is isomorphic to $\Gamma$? 

Cameron and Johnson \cite{CJ} give a sufficient condition for the random graph to arise as a Cayley graph of a countable group. We first show that $\mathbb{H}$ satisfies this condition and then apply their theorem to answer the existence question affirmatively. Indeed, if we include each inverse pair of nonidentity elements independently with probability $1/2$, the resulting Cayley graph is almost surely isomorphic to $\Gamma$. This abundance makes it natural to ask for a particular generating set specified by the group structure itself. We show that the set of all elements of order 7 works. Thus the Fra\"{i}ss\'e limit of finite graphs is a parameter-free, quantifier-free definable reduct of the Fra\"{i}ss\'e limit of finite groups.  


The proof returns to finite extensions. Given disjoint finite sets $U,V\subseteq\mathbb{H}$, put $F=\langle U\cup V\rangle$. Using small cancellation over free products and residual finiteness, we construct a finite overgroup of $F$ containing an element $x\notin F$ such that $a^{-1}x$ has order 7 for every $a\in U$, while $b^{-1}x$ has order different from 7 for every $b\in V$. Hall's extension property embeds this overgroup into $\mathbb{H}$ over $F$. The image of $x$ is then adjacent to every vertex in $U$ and to none in $V$, as required by the extension property of the random graph. To obtain the finite overgroup of $F$, we apply Martin and Steenbock's residual finiteness theorem \cite[Theorem~1.4]{MS} to a $C'(1/6)$
quotient of a free product of two finite groups. Their result combines
cubulation of the quotient with Agol's virtual specialness theorem
\cite{Agol}. 

For each integer $m\geq 2$, let $\Ha_m$ denote the set of elements of
$\Ha$ of order exactly $m$. With a few additional arguments, the same
proof works for every $m\geq 7$.

\begin{thmIntro}
For every integer $m\geq 7$, the Cayley graph $\Cay(\Ha,\Ha_m)$ is isomorphic to the random graph.
\end{thmIntro}

We also show that neither $\Cay(\Ha,\Ha_2)$ nor $\Cay(\Ha,\Ha_3)$
is isomorphic to the random graph. It remains open to determine
the least integer $m$ for which $\Cay(\Ha,\Ha_m)$ is isomorphic
to the random graph.

We use the order 7 Cayley graph $\Cay(\Ha,\Ha_7)$ to code finite
subsets of $\Ha$ uniformly. By passing to the cardinalities of the coded sets, we obtain a parameter-free interpretation of the standard model of arithmetic $\mathcal N=(\mathbb N,+,\cdot,0,1)$ in $\Ha$. This recovers Nies’s uniform coding of nonempty finite subsets \cite{Nies} and gives an alternative proof of Hodges’s interpretation theorem \cite{Hodges}.  Conversely, Hall's
construction \cite{Hal} gives a computable presentation of $\Ha$, which can be
interpreted in $\mathcal N$. The obstruction to a bi-interpretation
comes from their automorphisms. The structure $\mathcal N$ is rigid,
whereas every finite tuple in $\Ha$ is fixed pointwise by some
nontrivial inner automorphism.

The mutual interpretations above imply that $\Th(\Ha)$ admits no
axiomatization of bounded quantifier alternation. Consequently, there is
no fixed bound on the number of prenex quantifier blocks needed to
express all group formulas modulo $\Th(\Ha)$. In particular, $\Th(\Ha)$ is
not model complete. We show that $\Ha$ is existentially closed in every model of its
universal theory that contains it. Hence the failure of model
completeness yields the following.

\begin{thmIntro}
The first order theory of $\Ha$ does not admit a model companion.
\end{thmIntro}

Finally, we use the inner ultrahomogeneity of \(\Ha\) to obtain the following result.

\begin{thmIntro}
The first order theory of $\Ha$ does not geometrically eliminate imaginaries.
\end{thmIntro}

More generally, we show that if a group has infinitely many conjugacy classes and every nonidentity conjugacy class is infinite, then its pure group theory fails to geometrically eliminate imaginaries.

\section{Preliminaries}

\subsection{Hall's universal group}

Let $\Ha$ denote Hall's locally finite universal group. 

\begin{thm}[Hall {\cite{Hal}}]\label{Hall}
Hall's universal locally finite group $\Ha$ is the unique countable
locally finite group, up to isomorphism, with the following properties.
\begin{enumerate}
  \item Every finite group embeds in $\Ha$.
  \item Every isomorphism between finite subgroups of $\Ha$ extends to an
inner automorphism of $\Ha$.
\end{enumerate}
Moreover, for each integer $m\geq 2$, $\Ha=\Ha_m\Ha_m$, where $\Ha_m$ is the set of elements of order $m$.
\end{thm}

Hall's universal group can also be obtained as the Fra\"{i}ss\'e limit
of the class of all finite groups. In particular, it has the following
finite extension property.

\begin{cor}[Hall {\cite{Hal}}]\label{cor:Hall-extension}
Let $F\leq E$ be finite groups, and let $i\colon F\hookrightarrow\Ha$
be an embedding. Then there is an embedding
$\widehat{i}\colon E\hookrightarrow\Ha$ such that
$\widehat{i}|_F=i$.
\end{cor}

\begin{lemma}\label{lem:centerless}
Hall's universal group $\Ha$ is centerless.
\end{lemma}

\begin{proof}
Let $1\neq z\in\Ha$, and let $m$ be its order.
Since $C_m\times C_m$ embeds in $\Ha$, there is an
element $h\neq z$ of order $m$. By Theorem~\ref{Hall},
the isomorphism $\langle z\rangle\to\langle h\rangle$
sending $z$ to $h$ extends to an inner automorphism
of $\Ha$. Thus $z$ has a conjugate distinct from itself,
so $z\notin Z(\Ha)$.
\end{proof}

\subsection{Small cancellation quotients of free products}
Let
\[
P=A_1*A_2
\]
be the free product of two groups; $A_1$ and $A_2$ are called free factors. Recall that every nontrivial element
$g\in P$ has a unique reduced normal form
\[
g=s_1s_2\cdots s_n,
\]
where each syllable $s_i$ is a nonidentity element of one of the
free factors, and consecutive syllables belong to different free
factors. We call $n$ the syllable length of $g$, and write
$|g|=n$. A product $gh$ of elements of $P$ is called weakly reduced if
$g=1$ or $h=1$, or if the last syllable of $g$ and the first
syllable of $h$ are not mutual inverses. In particular, two
syllables at the join may combine into a single nonidentity
syllable.

\begin{defn}[Cyclically reduced and weakly cyclically
reduced elements]
Let $g\in P$ have free product normal form
\[
g=s_1\cdots s_n.
\]
If $n\leq 1$, then $g$ is both cyclically reduced and weakly
cyclically reduced. For $n\geq 2$, we say that $g$ is
\emph{cyclically reduced} if $s_1$ and $s_n$ lie in different
free factors, and \emph{weakly cyclically reduced} if
\[
s_ns_1\neq 1.
\]
\end{defn}

By the conjugacy theorem for free products
\cite[Chapter IV, Theorem 1.4]{LS},
two cyclically reduced elements of syllable length at least two
are conjugate in $P$ if and only if their normal forms are cyclic
permutations of one another.

A set $R$ of nontrivial weakly cyclically reduced elements of $P$
is called \emph{symmetrized} if it is closed under inversion and
under taking weakly cyclically reduced conjugates.

\begin{defn}[Piece]
Let $R$ be a symmetrized set of relators in $P$. A nontrivial
element $p\in P$ is called a \emph{piece} with respect to $R$ if
there are distinct $r_1,r_2\in R$ and $u_1,u_2\in P$
such that the products $r_1=pu_1$ and $r_2=pu_2$ are weakly reduced.
\end{defn}

\begin{defn}[The small cancellation condition
$C'(\lambda)$]
Let $\lambda>0$. A symmetrized set $R$ satisfies the
small cancellation condition $C'(\lambda)$ over $P$ if, for
every piece $p$ with respect to $R$ and every $r\in R$ that
admits a weakly reduced factorization $r=pu$, one has
\[
|p|<\lambda |r|.
\]
\end{defn}

We first recall two fundamental results from small cancellation
theory over free products.

\begin{thm}[{\cite[Chapter V, Corollary 9.4]{LS}}]
\label{thm:free-factor-embedding}
Let $P=A_1*\cdots*A_m$, and let $R$ be a symmetrized set of
relators satisfying $C'(1/6)$ over $P$ such that $|r|>6$ for every $r\in R$. Then the natural
projection
\[
P\longrightarrow P/\langle\!\langle R\rangle\!\rangle
\]
is injective on each free factor $A_i$, where $\langle\!\langle R\rangle\!\rangle$
 is the normal closure of $R$.
\end{thm}

We shall also use the following consequence of Greendlinger's
lemma for free products
\cite[Chapter V, Theorem 9.3]{LS}.

\begin{cor}\label{thm:greendlinger-consequence}
Let $R$ be a nonempty symmetrized set of relators in $P$
satisfying $C'(1/6)$, and let $L=\min\{|r|:r\in R\}$.
Suppose that $L>6$. Then every nontrivial element of $P$
of syllable length at most $L/2$ has nontrivial image in $P/\langle\!\langle R\rangle\!\rangle$.
\end{cor}

\begin{defn}[Residual finiteness]
A group $G$ is called \emph{residually finite} if, for every
nonidentity $g\in G$, there exists a homomorphism
from $G$ to a finite group such that the image of $g$ is nonidentity.
\end{defn}

We next record a result that establishes residual finiteness
for certain small cancellation quotients of free products.

\begin{thm}[Martin and Steenbock, {\cite[Theorem 1.4]{MS}}]
\label{thm:MS-residual-finiteness}
Let
\[
P=G_1*\cdots*G_m,
\]
where each $G_i$ is hyperbolic and cubulable. Let $R\subseteq P$
be finite, and let $\overline{R}$ be its symmetrization.
Suppose that $\overline{R}$ satisfies $C'(1/6)$ over $P$ and
that every $r\in\overline{R}$ has syllable length at least $6$.
Then $Q=P/\langle\!\langle\overline{R}\rangle\!\rangle$ is residually finite and linear over $\mathbb Z$, and every
quasiconvex subgroup of $Q$ is separable.
\end{thm}

\begin{rmk}
Martin and Steenbock's main theorem \cite[Theorem~1.1]{MS}
gives a proper, cocompact action of $Q$ on a
$\mathrm{CAT}(0)$ cube complex. The hyperbolicity of the
free factors also ensures that $Q$ is word-hyperbolic
\cite[Theorem~2.2]{MS}. They then invoke Agol's celebrated
virtual specialness theorem \cite[Theorem~1.1]{Agol}. It follows that $Q$ has a finite-index subgroup that embeds in a right-angled Artin group. Since right-angled Artin groups are residually finite, so is $Q$.
\end{rmk}

We will only need the residual finiteness
statement of Theorem~\ref{thm:MS-residual-finiteness}.

\begin{cor}
\label{cor:finite-factors-residually-finite}
Let
\[
P=G_1*\cdots*G_m
\]
be a free product of finite groups. Let $R\subseteq P$ be
finite, and let $\overline{R}$ be its symmetrization. Suppose
that $\overline{R}$ satisfies $C'(1/6)$ over $P$ and that every
$r\in\overline{R}$ has syllable length at least $6$. Then $Q=P/\langle\!\langle\overline{R}\rangle\!\rangle$ is residually finite.
\end{cor}

\begin{proof}
Every finite group is word-hyperbolic and admits a proper,
cocompact action on a point, which is a CAT(0) cube complex.
Thus every $G_i$ is hyperbolic and cubulable. The conclusion
follows immediately from
Theorem~\ref{thm:MS-residual-finiteness}.
\end{proof}

\section{Cayley graphs of Hall's universal group}

In this section we will prove the main theorem. We will first show, using a condition of Cameron and Johnson, that Hall's universal group admits a Cayley graph isomorphic to the random graph. In fact, if $S\subseteq\Ha\setminus\{1\}$ is chosen by including each inversion orbit independently with probability $1/2$,
then $\Cay(\Ha,S)$ is almost surely isomorphic to the random graph. We then show that, for every integer $m\geq 7$, taking the elements of order exactly $m$ as the connection set yields a Cayley graph isomorphic to the random graph.

For a group $G$ and an element $g\in G$, the \emph{square root set}
of $g$ is
\[
    \sqrt{g}=\{x\in G\mid x^2=g\}.
\]
This set is called \emph{principal} if $g=1$ and
\emph{non-principal} otherwise.

\begin{defn}[Cameron and Johnson {\cite{CJ}}]
A group $G$ has \emph{property $P'$} if the complement of every finite
union of translates of non-principal square root sets is infinite.
\end{defn}

\begin{thm}[Cameron and Johnson {\cite[Theorem~3.1]{CJ}}]\label{CJRandom}
Let $G$ be a countable group with property $P'$. Then $\Aut(\Gamma)$
contains a regular subgroup isomorphic to $G$.
\end{thm}

Identifying the vertices of $\Gamma$ with $G$ via the orbit map of
a fixed vertex gives a Cayley graph of $G$ isomorphic to $\Gamma$.

In fact, property $P'$ gives the stronger conclusion that a randomly
chosen symmetric subset $S\subseteq G\setminus\{1\}$ almost surely
yields a Cayley graph $\operatorname{Cay}(G,S)$ isomorphic to $\Gamma$.

\begin{thm}[{\cite[Proposition~8.3]{Cam}}]
Let $G$ be a countable group with property $P'$. Then there is an
inverse-closed generating set $S\subseteq G\setminus\{1\}$ such that
$\operatorname{Cay}(G,S)$ is isomorphic to $\Gamma$. Moreover, if $S$
is chosen by including each inversion orbit in $G\setminus\{1\}$
independently with probability $1/2$, then
$\operatorname{Cay}(G,S)$ is almost surely isomorphic to $\Gamma$.
\end{thm}

We now verify that Hall’s group has property \(P'\).

\begin{thm}
Hall's universal group can be embedded as a regular subgroup of $\Aut(\Gamma)$, and thus the random graph can be constructed as a Cayley graph of Hall’s universal group.
\end{thm}
\begin{proof}
To verify $P'$, let $a_1,\ldots,a_n,b_1,\ldots,b_n\in\Ha$ and let
$D\subseteq\Ha$ be finite. We will find an element outside
$D\cup\bigcup_{i=1}^n a_i\sqrt{b_i}$, where
$\sqrt{b_i}=\{y\in\Ha:y^2=b_i\}$.

By local finiteness, there is a finite subgroup $F\leq\Ha$ containing
$D$ and all the $a_i$ and $b_i$. By the finite extension property of
$\Ha$, the group $F\times C_3$ embeds in $\Ha$ over $F$. Let $x$ be the
image of a generator of the second factor. Then $x\notin F$, so
$x\notin D$. Moreover, for each $i$,
\[
(a_i^{-1}x)^2=a_i^{-2}x^2\notin F,
\]
since $x$ commutes with $F$ and $x^2$ is a nonidentity element of the
second factor. As $b_i\in F$, we have $(a_i^{-1}x)^2\ne b_i$, so
$x\notin a_i\sqrt{b_i}$. Since $D$ was arbitrary, the complement of
$\bigcup_{i=1}^n a_i\sqrt{b_i}$ is infinite, as required.
\end{proof}

A randomly chosen symmetric subset of $\Ha\setminus\{1\}$ almost
surely yields a Cayley graph isomorphic to the random graph, but
finding an explicit connection set is less straightforward. The
sets $\Ha_m$ are natural candidates. For each $m\geq 2$, the set
$\Ha_m$ is closed under inversion and generates $\Ha$ by
Theorem~\ref{Hall}. The first two choices both fail the extension
property of the random graph, so neither $\Cay(\Ha,\Ha_2)$ nor
$\Cay(\Ha,\Ha_3)$ is isomorphic to it.

\begin{fact}
The Cayley graph $\Cay(\Ha,\Ha_2)$ is not isomorphic to the random graph.
\end{fact}
\begin{proof}
Let $g$ have order $3$. If the graph satisfied the extension property,
there would be a vertex $x\notin\{1,g,g^2\}$ adjacent to $1$ and $g$ but
not to $g^2$. Then $x^2=(gx^{-1})^2=1$, so $xgx=g^{-1}$. Hence
\[
(g^2x^{-1})^2=g^2xg^2x=g^2g^{-2}=1.
\]
Since $x\ne g^2$, the element $g^2x^{-1}$ has order $2$. Thus $x$ is
adjacent to $g^2$, a contradiction.
\end{proof}

\begin{fact}
The Cayley graph $\Cay(\Ha,\Ha_3)$ is not isomorphic to the random graph.
\end{fact}
\begin{proof}
Let $g$ have order $4$. If the graph satisfied the extension property,
there would be a vertex $y\notin\{1,g,g^2,g^3\}$ adjacent to $1$, $g$,
and $g^2$ but not to $g^3$. Put $h=y^{-1}$ and $a_i=h^igh^{-i}$ for
$i=0,1,2$. Then $h^3=(gh)^3=(g^2h)^3=1$. Since
\[
(g^j h)^3=a_0^j a_1^j a_2^j
\]
for every integer $j$, taking $j=1,2$ gives
$a_0a_1=a_2^{-1}$ and $a_0^2a_1^2=a_2^{-2}=(a_0a_1)^2$.
Thus $a_0$ and $a_1$ commute, and so does $a_2=(a_0a_1)^{-1}$ with both.
Taking $j=-1$ now gives
\[
(g^3y^{-1})^3=(g^{-1}h)^3
=a_0^{-1}a_1^{-1}a_2^{-1}=(a_0a_1a_2)^{-1}=1.
\]
Since $y\ne g^3$, the element $g^3y^{-1}$ has order $3$. Thus $y$ is
adjacent to $g^3$, a contradiction.
\end{proof}

\begin{thm}
The Cayley graph $Cay(\Ha, \Ha_7)$ is isomorphic to the random graph.   
\end{thm}

\begin{proof}
Let $U,V$ be disjoint finite subsets of $\Ha$, and put
$F=\langle U\cup V\rangle$. Since $\Ha$ is locally finite, $F$ is
finite. We construct a finite overgroup $E$ of $F$ and an element
$x\in E\setminus F$ such that
\[
  \operatorname{ord}(u^{-1}x)=7 \quad (u\in U),
  \qquad
  (v^{-1}x)^8=1 \quad (v\in V).
\]
Let $c,d,t$ generate cyclic groups of orders $2,3,7$, respectively.
Set $K=F\times\langle c\rangle$ and put
\[
  P=(K\times\langle d\rangle)*\langle t\rangle,
  \qquad w=dt(d^2t)^3.
\]
Let $R$ denote the set of relators $\{(u^{-1}w)^7~(u\in U), (v^{-1}w)^8~(v\in V), (cw)^7\},$ and let $\overline{R}$ be its symmetrization.
Form the quotient
\[
Q=P/\langle\!\langle\overline{R}\rangle\!\rangle.
\]

We verify that $\overline{R}$ satisfies
$C'(1/6)$ over $P$. For each $a\in U^{-1}\cup V^{-1}\cup\{c\}$, the word $aw$ is cyclically
reduced and has syllable sequence
\[
  ad,\ t,\ d^2,\ t,\ d^2,\ t,\ d^2,\ t.
\]
The defining relators therefore have syllable length $56$ or $64$.
Any eight consecutive syllables in a cyclic reading of a defining
relator or its inverse determine the orientation, since the
syllables from $\langle t\rangle$ are all $t$ in the positive
orientation and all $t^{-1}$ in the negative orientation. The
same eight syllables contain exactly one occurrence of $ad$ or
$(ad)^{-1}$, according to the orientation. Its
$\langle d\rangle$-coordinate distinguishes it from the other
syllables in the first factor. Its value determines $a$, and its
location determines the position modulo 8. Since $a$ also determines the exponent of the defining relator,
the full cyclic reading is uniquely specified.

Every relator in $\overline{R}$ has the form
\[
  y s_2\cdots s_L z,
\]
where $s_1\cdots s_L$ is a cyclic permutation of a defining relator
or of its inverse, $y$ and $z$ lie in the factor containing
$s_1$ with $zy=s_1$ and $y\neq 1$. Here $z$ is omitted if it is identity. Indeed, if the first
and last syllables of a weakly cyclically reduced conjugate lie in
the same factor, moving the last syllable to the front and combining
them gives a cyclically reduced conjugate. By the free product
conjugacy theorem, this is a cyclic permutation of the defining
relator or its inverse. In particular, every relator in $\overline{R}$ has
length at least $56$.

Suppose that a piece $p$ has length at least 10. Then there are two distinct 
$r_1, r_2\in\overline{R}$ such that $r_1=pu_1$ and $r_2=pu_2$, with both products weakly reduced. Thus, the first nine syllables agree, since only the last syllable
of $p$ can combine with the next syllable. The common syllables in
positions 2 through 9 are therefore eight unchanged syllables
in the two cyclic readings. By the preceding calculation, they
determine the same cyclic permutation $s_1\cdots s_L$. The initial
syllable $y$ also agrees, so $z=s_1y^{-1}$ is the same in both
relators. This contradicts their distinctness. Hence every piece
has length at most $9<56/6$, satisfying $C'(1/6)$.

By Theorem~\ref{thm:free-factor-embedding}, the first free factor embeds in $Q$.
Both free factors are finite, so
Corollary~\ref{cor:finite-factors-residually-finite} gives residual finiteness
of $Q$. Taking the image of the product of finitely many finite
quotients that separate the nonidentity elements of $K$ from the
identity, we obtain a finite quotient $q\colon Q\twoheadrightarrow E$
that is injective on $K$. Identify $K$ with its image in $E$, and
let $x$ be the image of $w$.

The relation $(cx)^7=1$ forces $x\notin F$. Otherwise, suppose that $x\in F$. Since $K=F\times\langle c\rangle$, the element $(cx)^7=c^7x^7=cx^7\neq 1$.
For each $u\in U$, since $x\notin F$, the element $u^{-1}x\neq 1$. Then the relation $(u^{-1}x)^7=1$ implies that it has order exactly $7$. For each $v\in V$, the relation
$(v^{-1}x)^8=1$ implies that the order of $v^{-1}x$ divides $8$
and is therefore different from $7$.

By Hall's finite extension property, $E$ embeds into $\Ha$ by an
embedding fixing $F$ pointwise. The image of $x$ lies outside
$U\cup V$ and is adjacent in $\Cay(\Ha,\Ha_7)$ to every vertex
of $U$ and to no vertex of $V$. Thus $\Cay(\Ha,\Ha_7)$ satisfies
the extension property of the random graph.
\end{proof}

The preceding construction extends to every integer $m\geq 7$. Choose a prime $q>m$ and take the auxiliary element $c$ to have order $q$ in $K=F\times\langle c\rangle$. Keeping the same $d,t,w$, impose the relations
\[
(u^{-1}w)^m=1\quad(u\in U),\qquad
(v^{-1}w)^q=1\quad(v\in V),\qquad
(cw)^m=1.
\]
The same piece estimate applies. Every piece has syllable length at most $9$, whereas every defining relator has syllable length at least $8m$. Since $9<8m/6$, the relators satisfy $C'(1/6)$. As before, $K$ embeds in the resulting quotient $Q$, and $Q$ is residually finite.

Every relator in the symmetrized set has syllable length at least
$8m$. For each $u\in U$ and each prime divisor $\ell$ of $m$,
the word $(u^{-1}w)^{m/\ell}$ has syllable length $8m/\ell\leq 4m$.
By Corollary~\ref{thm:greendlinger-consequence}, it remains nontrivial in $Q$. By residual finiteness, we can take a finite quotient of $Q$ in which $K$ embeds and all these words remain nontrivial. We identify $K$ with its image and let $x$ be the image of $w$. For every $u\in U$, the defining relation gives $(u^{-1}x)^m=1$. On the other hand, for every prime divisor $\ell$ of $m$, the element $(u^{-1}x)^{m/\ell}\neq 1$. Thus, the order of $(u^{-1}x)$ is exactly $m$.

If $x\in F$, then $cx\in K$, and its image in $K/F$ would have order $q$, contradicting $(cx)^m=1$. Thus $x\notin F$. In particular, for every $v\in V$, the element $v^{-1}x$ is nontrivial and has order dividing $q$, so it has order exactly $q$. Embedding the finite quotient into $\Ha$ over $F$ gives the required witness to the extension property.

We have therefore proved the following.

\begin{thm}
For every integer $m\geq 7$, let $\Ha_m$ denote the set of elements of $\Ha$ of order $m$. Then $\Cay(\Ha,\Ha_m)$ is isomorphic to the random graph. 
\end{thm}

\section{Interpretability}
Hall's construction \cite{Hal} gives a solvable word problem for $\Ha$, so
$\Ha$ admits a computable presentation. It follows that $\Ha$ is
interpretable without parameters in $(\mathbb{N},+,\cdot,0,1)$.

To interpret standard arithmetic in $\Ha$, we use the following lemma, which follows from Hodges \cite[Lemma~8]{Hodges}. Rzepecki \cite[Proposition~4.33]{Rzepecki} proves a more general version. We include a proof along the lines of Macintyre's argument for existentially closed groups \cite[Lemma~4]{Mac}.

\begin{lemma}\label{UniformDef}
For each $n\geq 1$, there is a universal group formula
$\varphi_n(x_1,\ldots,x_n,y)$ such that, for all
$a_1,\ldots,a_n,b\in\Ha$,
\[
  \Ha\models\varphi_n(a_1,\ldots,a_n,b)
  \quad\Longleftrightarrow\quad
  b\in\langle a_1,\ldots,a_n\rangle.
\]
\end{lemma}

\begin{proof}
Put $F=\langle a_1,\ldots,a_n\rangle$ and define
\[
  \varphi_n(x_1,\ldots,x_n,y)
  \;=\;
  \forall z\left(
    \bigwedge_{i=1}^{n} x_i^z=x_i \ \longrightarrow\ y^z=y
  \right).
\]
Every element centralizing $a_1,\ldots,a_n$ centralizes $F$
pointwise. Thus $\Ha\models\varphi_n(a_1,\ldots,a_n,b)$
whenever $b\in F$.

Suppose now that $b\notin F$, and put $G=\langle F,b\rangle$.
By local finiteness, $G$ is finite. Disjoint amalgamation for
finite groups gives a finite overgroup $E$ of $G$ and an embedding
$j\colon G\hookrightarrow E$ fixing $F$ pointwise such that
$G\cap j(G)=F$. By Hall's finite extension property, we may regard
$E$ as a subgroup of $\Ha$ containing the original copy of $G$.

By inner ultrahomogeneity, the isomorphism $j\colon G\to j(G)$ is
induced by conjugation by some $c\in\Ha$. Since $j$ fixes $F$
pointwise, we have $a_i^c=a_i$ for every $i$. On the other hand,
$b^c=j(b)\ne b$, since $b\notin F$ and $G\cap j(G)=F$. Hence
$\Ha\not\models\varphi_n(a_1,\ldots,a_n,b)$.
\end{proof}

The proof shows that the double centralizer of any finite tuple $\bar a$
in $\Ha$ is the subgroup generated by $\bar a$:
\[
C_{\Ha}\bigl(C_{\Ha}(\bar a)\bigr)=\langle\bar a\rangle.
\]

Hodges \cite[Theorem~7(a)]{Hodges} interpreted standard arithmetic in $\Ha$, using element orders to represent natural numbers and group-word interpolation to code finite sequences. Using Hodges's results, Nies \cite[Lemma~4.2]{Nies} gave a uniform coding of the nonempty finite subsets of $\Ha$. We give an alternative proof of Hodges's theorem by using the random graph reduct to code all finite subsets of $\Ha$ and representing natural numbers by their cardinalities.

\begin{thm}
Standard arithmetic $(\mathbb{N},+,\cdot,0,1)$ is
interpretable without parameters in Hall's universal locally finite
group $\Ha$.
\end{thm}

\begin{proof}
We represent natural numbers by finite subsets of $\Ha$ modulo
equicardinality.
By the preceding lemma, $\varphi_2(d_1,d_2,x)$ defines membership in
$\langle d_1,d_2\rangle$, uniformly in $d_1$ and $d_2$.
Every finite subset $A$ of $\Ha$ is contained in a 2-generated
finite subgroup. Indeed, $F=\langle A\rangle$ is finite, and its
regular action embeds it into $\operatorname{Sym}(F)$, which is
generated by at most two elements. By Hall's finite extension
property, this symmetric group embeds into $\Ha$ over $F$.

Let $\sim$ denote adjacency in $\Cay(\Ha,\Ha_7)$. A triple
$\bar d=(d_1,d_2,z)$ codes the set
\[
  S_{\bar d}=\{x\in\Ha:\varphi_2(d_1,d_2,x)\textrm{ and } x\sim z\}.
\]
This set is finite. Conversely, given a finite subset $A$ of $\Ha$,
choose $D=\langle d_1,d_2\rangle$ containing $A$. The extension
property of the random graph gives a vertex $z$ adjacent to every
element of $A$ and to none of $D\setminus A$. Then
$S_{(d_1,d_2,z)}=A$. Thus quantification over finite subsets can
be expressed by quantification over triples, with membership given
by the formula above.

To define equicardinality, let $A,B$ be finite subsets of $\Ha$
and put $F=\langle A\cup B\rangle$. Embed $F\times F$ into $\Ha$
over its first factor. Since the two factors are conjugate in
$\Ha$, there is a $g\in\Ha$ such that $F^g$ is the second factor,
where $b^g=g^{-1}bg$. In particular, the map
\[
  A\times B\longrightarrow\Ha,\qquad (a,b)\longmapsto ab^g
\]
is injective. Thus finite subsets of $A\times B$ can be
encoded by finite subsets of $\Ha$.

Let $\bar p$ and $\bar q$ be codes for $A$ and $B$, respectively,
and write $\mu(x;\bar r)$ for membership in the finite set
coded by $\bar r$. Define $A\equiv_{\#}B$ by requiring that
\[
\begin{aligned}
\Ha\models \exists g\,\exists\bar r\;\bigl[&
\forall a\bigl(\mu(a;\bar p)\rightarrow
  \exists!b\,(\mu(b;\bar q)\land\mu(ab^g;\bar r))\bigr)
\\
&{}\land
\forall b\bigl(\mu(b;\bar q)\rightarrow
  \exists!a\,(\mu(a;\bar p)\land\mu(ab^g;\bar r))\bigr)
\bigr].
\end{aligned}
\]
The first clause
assigns to each element of $A$ exactly one element of $B$,
and the second ensures that each element of $B$ has exactly
one preimage. Thus any witnesses $g,\bar r$ define a bijection
from $A$ to $B$.

Conversely, suppose that $f:A\to B$ is a bijection. Choose $g$
as above so that the map $(a,b)\mapsto ab^g$ is injective on
$A\times B$, and let $\bar r$ code the finite set
\[
  P=\{af(a)^g:a\in A\}.
\]
Then, for $a\in A$ and $b\in B$, we have $ab^g\in P$ if and
only if $b=f(a)$, so the two clauses hold. Consequently,
$A\equiv_{\#}B$ if and only if $|A|=|B|$. In particular,
equicardinality is first order definable on the codes.
Write $\operatorname{Eq}(\bar p,\bar q)$ for the formula defining
equicardinality on codes, and let $\bar p,\bar q,\bar r$ code
$A,B,C$, respectively.

For addition, define
\[
\begin{aligned}
\operatorname{Add}(\bar p,\bar q,\bar r)
\ :=\ \exists\bar u\,\exists\bar v\;\bigl[&
\operatorname{Eq}(\bar u,\bar p)
\land\operatorname{Eq}(\bar v,\bar q)
\\
&{}\land\forall x\bigl(\mu(x;\bar r)\leftrightarrow
  (\mu(x;\bar u)\lor\mu(x;\bar v))\bigr)
\\
&{}\land\forall x\,\neg
  (\mu(x;\bar u)\land\mu(x;\bar v))\bigr].
\end{aligned}
\]
Thus the sets coded by $\bar u,\bar v$ partition $C$ and have
cardinalities $|A|,|B|$, respectively. Conversely, if $|C|=|A|+|B|$, choose a subset $C_0\subseteq C$
of cardinality $|A|$ and put $C_1=C\setminus C_0$. Then
$|C_1|=|B|$. Since every finite subset of $\Ha$ has a code,
we may choose codes $\bar u,\bar v$ for $C_0,C_1$, respectively,
which witness $\operatorname{Add}(\bar p,\bar q,\bar r)$.
Hence
\[
  \Ha\models\operatorname{Add}(\bar p,\bar q,\bar r)
  \quad\Longleftrightarrow\quad
  |C|=|A|+|B|.
\]

For multiplication, first express injectivity of the map
$(a,b)\mapsto ab^g$ on $A\times B$ by
\[
\begin{aligned}
\operatorname{Inj}(\bar p,\bar q,g)
\ :=\ \forall a\,\forall a'\,\forall b\,\forall b'\;\bigl[
&\bigl(\mu(a;\bar p)\land\mu(a';\bar p)
\\
&{}\land\mu(b;\bar q)\land\mu(b';\bar q)
\\
&{}\land ab^g=a'(b')^g\bigr)
\rightarrow(a=a'\land b=b')\bigr].
\end{aligned}
\]
Now define
\[
\begin{aligned}
\operatorname{Mult}(\bar p,\bar q,\bar r)
\ :=\ \exists g\,\exists\bar s\;\bigl[&
\operatorname{Inj}(\bar p,\bar q,g)
\land\operatorname{Eq}(\bar s,\bar r)
\\
&{}\land\forall x\bigl(\mu(x;\bar s)\leftrightarrow
  \exists a\,\exists b\,
  (\mu(a;\bar p)\land\mu(b;\bar q)\land x=ab^g)\bigr)
\bigr].
\end{aligned}
\]
The set coded by $\bar s$ is exactly the image of $A\times B$
under the injective map $(a,b)\mapsto ab^g$, so it has cardinality
$|A|\cdot|B|$. Conversely, the preceding construction supplies
a suitable $g$, and this finite image has a code. Therefore
\[
  \Ha\models\operatorname{Mult}(\bar p,\bar q,\bar r)
  \quad\Longleftrightarrow\quad
  |C|=|A|\cdot|B|.
\]

Finally, the classes of the empty set and of the singletons
define $0$ and $1$, respectively. Every finite cardinality occurs,
so the map sending the class of a code $\bar d$ to $|S_{\bar d}|$
is an isomorphism from the interpreted structure to
$(\mathbb{N},+,\cdot,0,1)$. All the defining formulas are
parameter-free.
\end{proof}

Although both structures are mutually interpretable, there is no bi-interpretation even if we allow parameters. Note that a parameter-free bi-interpretation induces an isomorphism of automorphism groups (see \cite[p.~125]{AKNS}). If $M$ is bi-interpretable, even with parameters, with a rigid
structure $N$ (that is, $\operatorname{Aut}(N)=\{1\}$), then
$\operatorname{Aut}(M/\bar a)=\{1\}$ for some finite tuple
$\bar a$ in $M$.

Indeed, let $N^*$ be the copy of $N$ interpreted in $M$, and
let $M^{**}$ be the copy of $M$ obtained by interpreting $M$
in $N^*$. Bi-interpretability provides a definable isomorphism
\[
  h:M\longrightarrow M^{**}.
\]
Choose a finite tuple $\bar a$ in $M$ containing the parameters
needed for these interpretations and for $h$, including
representatives of any parameters from $N^*$.

Let $\sigma\in\operatorname{Aut}(M/\bar a)$. It induces an
automorphism $\sigma_N$ of $N^*$, which is the identity since
$N^*\cong N$ is rigid. By construction of the composite
interpretation, the action $\sigma_{**}$ of $\sigma$ on
$M^{**}$ coincides with the action induced by $\sigma_N$.
Consequently, $\sigma_{**}$ is also the identity.

Since $h$ is $\bar a$-definable, for every $m\in M$ we have
\[
  h(\sigma(m))=\sigma_{**}(h(m))=h(m).
\]
Injectivity of $h$ gives $\sigma(m)=m$. Thus
$\operatorname{Aut}(M/\bar a)=\{1\}$.

\begin{thm}
Hall's universal locally finite group $\Ha$ and the standard model
of arithmetic $\mathcal N$ are not bi-interpretable, even with
parameters.
\end{thm}

\begin{proof}
Every element of $\mathcal N$ is denoted by a numeral, so
$\Aut(\mathcal N)=\{1\}$. By the discussion in the preceding paragraph,
a bi-interpretation would therefore make $\Aut(\Ha/\bar a)$ trivial
for some finite tuple $\bar a$ in $\Ha$.

Let $\bar a$ be any finite tuple in $\Ha$ and put
$F=\langle\bar a\rangle$. By local finiteness, $F$ is finite.
Hall's finite extension property gives an embedding of
$F\times S_2$ into $\Ha$ fixing the first factor pointwise.
Let $s$ be the image of the nonidentity element of the second
factor. Then $s$ centralizes $F$, so conjugation by $s$ fixes
$\bar a$ pointwise. Since $\Ha$ is centerless  and $s\ne 1$, this
conjugation is a nontrivial automorphism. Thus
$\Aut(\Ha/\bar a)\ne\{1\}$ for every finite tuple $\bar a$,
a contradiction.
\end{proof}

\section{Quantifier complexity and model companions}

By Hodges's theorem, $\Ha$ and standard arithmetic are mutually
interpretable without parameters. Since the language of groups is finite,
we have that $\Th(\Ha)$ admits no axiomatization of
bounded quantifier alternation. For each $k$, let $B_k$
be the set of prenex group sentences with at most $k$
quantifier blocks that hold in $\Ha$. Since $\Ha$ has a
computable presentation, $B_k$ is arithmetically definable.
If $\Th(\Ha)$ admitted an axiomatization of bounded quantifier alternation, then some $B_k$ would axiomatize $\Th(\Ha)$.
Consequently, provability from $B_k$ would make $\Th(\Ha)$ arithmetically definable. However, the interpretation of standard arithmetic in $\Ha$ gives a computable reduction of true arithmetic to $\Th(\Ha)$, contradicting Tarski's undefinability theorem. Thus $\Th(\Ha)$ admits no axiomatization
of bounded quantifier alternation.
Consequently, there is no fixed $k$ such
that every group formula is equivalent modulo $\Th(\Ha)$ to a prenex
formula with at most $k$ quantifier blocks. In particular, $\Th(\Ha)$ is
not model complete.

The next lemma is an auxiliary result towards showing that $\Th(\Ha)$ does not admit a model companion.

\begin{lemma}
Let $U=\Th_{\forall}(\Ha)$. Then $\Ha$ is an existentially closed model of $U$.
\end{lemma}

\begin{proof}
Suppose $\Ha\leq G\models U$ and
$G\models\exists\bar y\,\theta(\bar a,\bar y)$, where $\bar a$ is a tuple
from $\Ha$ and $\theta$ is quantifier-free. By local finiteness,
$F=\langle\bar a\rangle$ is finite. Consider the existential sentence
asserting that there is a copy of the multiplication table of $F$
together with a witness to $\theta$ for the tuple corresponding to
$\bar a$. This sentence holds in $G$, hence also in $\Ha$, since otherwise
its universal negation would belong to $U$.

The obtained copy of $F$ is isomorphic to the original $F$ by the map
matching their enumerations. Inner ultrahomogeneity extends this
isomorphism to an automorphism of $\Ha$, which carries the witness to
one over the original tuple $\bar a$.
\end{proof}

\begin{cor}
The first order theory of $\Ha$ does not admit a model companion.    
\end{cor}
\begin{proof}
Suppose that $\Th(\Ha)$ has a model companion $T^*$. Then $T^*$ is
also a model companion of $U$, so its models are exactly the existentially
closed models of $U$. By the preceding lemma, $\Ha\models T^*$. Now since $T^*$ is model complete we get that $\Th(\Ha)$ is model complete.
\end{proof}

\section{Imaginaries}

Although true arithmetic eliminates imaginaries, $\Th(\Ha)$ fails even
geometric elimination of imaginaries. 

Recall that a complete first order theory $T$ has geometric elimination
of imaginaries if, in a monster model of $T$, every imaginary element
$e$ is interalgebraic with a finite tuple $\bar a$ of real elements,
that is,
\[
e\in acl^{\mathrm{eq}}(\bar a)
\qquad\text{and}\qquad
\bar a\in acl^{\mathrm{eq}}(e).
\]

\begin{thm}\label{thm:conjugacy-imaginaries}
Let $G$ be a group with infinitely many conjugacy classes.
Suppose that every nonidentity element of $G$ has an infinite
conjugacy class. Then $\operatorname{Th}(G)$, in the pure group
language, does not geometrically eliminate imaginaries.
\end{thm}

\begin{proof}
Let $\mathfrak C$ be a monster model of $\operatorname{Th}(G)$.
Both hypotheses are expressed by first order schemes and
therefore hold in $\mathfrak C$.
Consider the definable equivalence relation
\[
    E(x,y) := \exists g\;(y=g^{-1}xg).
\]
Since $\mathfrak C/E$ is infinite, saturation gives
\[
    e\in(\mathfrak C/E)
    \setminus\operatorname{acl}^{\mathrm{eq}}(\varnothing).
\]
Every inner automorphism fixes every conjugacy class setwise,
so
\[
    \operatorname{Inn}(\mathfrak C)
    \leq \operatorname{Aut}(\mathfrak C/e).
\]
If $b\in\operatorname{acl}^{\mathrm{eq}}(e)\cap\mathfrak C$,
then its orbit under $\operatorname{Aut}(\mathfrak C/e)$ is
finite. Its orbit under inner automorphisms, which is its
conjugacy class, is therefore finite as well. By hypothesis,
$b=1$. Hence
\[
    \operatorname{acl}^{\mathrm{eq}}(e)\cap\mathfrak C=\{1\}.
\]

If a finite real tuple $\bar a$ were interalgebraic with $e$,
every coordinate of $\bar a$ would consequently equal $1$.
Since the identity is parameter-free definable, this would give
\[
    e\in\operatorname{acl}^{\mathrm{eq}}(\bar a)
    =\operatorname{acl}^{\mathrm{eq}}(\varnothing),
\]
a contradiction.
\end{proof}

\begin{cor}
The theory of Hall's universal locally finite group $\Ha$
does not geometrically eliminate imaginaries.
\end{cor}

\begin{proof}
Since every finite group embeds in $\Ha$, the group $\Ha$
contains elements of every finite order. Elements of different
orders are not conjugate, so $\Ha$ has infinitely many
conjugacy classes.
Now let $a\in\Ha\setminus\{1\}$, and let $n$ be its order.
For every $r\geq 1$, the finite group
$(\mathbb Z/n\mathbb Z)^r$ embeds in $\Ha$, giving at least
$r$ distinct elements of order $n$. Every such element is
conjugate to $a$. Indeed, the isomorphism between their cyclic
subgroups sending one generator to the other extends to an
inner automorphism of $\Ha$. Thus the conjugacy class of $a$
is infinite.
The conclusion follows from
Theorem~\ref{thm:conjugacy-imaginaries}.
\end{proof}

\subsection*{Acknowledgement}
The second named author is grateful to Peter Cameron for bringing the reference \cite{CJ} to his attention at the Third International Conference on Group Actions and Transitive Graphs in 2018, and to Andr\'{e} Nies for pointing out the reference \cite{Hodges} at the Tianyuan Workshop on Computability Theory and Descriptive Set Theory in 2025. 

\subsection*{Declaration of generative AI and AI-assisted technologies in the manuscript preparation process}

We used OpenAI's ChatGPT to generate and explore mathematical ideas
and to help draft and revise the exposition. In particular, we developed
the proof of Theorem~3.7 through a series of exchanges with ChatGPT.
We reviewed the tool's output, checked the mathematical arguments,
and wrote the proofs presented in this paper. We take full responsibility
for the content of the manuscript.

\vspace{2cm}

\noindent
\begin{minipage}[t]{0.6\textwidth}
\raggedright
\textbf{Rizos Sklinos}\\[1em]
Academy of Mathematics \& Systems Science\\
Chinese Academy of Sciences\\
Room 1016, SiYuan Building\\
No.~55 Zhongguancun East Road\\
Haidian District, Beijing 100190\\
P.~R.~China
\end{minipage}\hfill
\begin{minipage}[t]{0.4\textwidth}
\raggedright
\textbf{Shichang Song}\\[1em]
School of Mathematics and Statistics\\
Beijing Jiaotong University\\
No.~3 Shangyuancun\\
Haidian District, Beijing 100044\\
P.~R.~China
\end{minipage}

\end{document}